\documentstyle[12pt]{amsart}

\newtheorem{th}{Theorem}[section]

\newcommand{\rar}{\rightarrow} 	
	
\newcommand{\bfp}{{\Bbb{P}}}   	\newcommand{\bfc}{{\Bbb{C}}}	
\newcommand{\bfr}{{\Bbb{R}}}	\newcommand{\bfh}{{\Bbb{H}}}
\newcommand{\bfq}{{\Bbb{Q}}}	
\newcommand{\bfz}{{\Bbb{Z}}}

\newcommand{\sym}{{\operatorname{Sym}}}

\newcommand{\Ga}{\Gamma}	
\newcommand{\GH}{\Gamma\setminus \bfh}

\begin{document}

\noindent%
{\large  \bf
 A linear lower bound on the gonality of modular curves \\[2mm]
}

\noindent%
{Dan Abramovich\footnote{Partially supported by NSF grant
DMS-9503276 and an Alfred P. Sloan research Fellowship.} \\ 
Department of Mathematics, Boston University\\
111 Cummington, Boston, MA 02215, USA \\
{\tt abrmovic@@math.bu.edu}}\\[2mm]
Preliminary version, \today

\addtocounter{section}{-1}
\section{INTRODUCTION}
\subsection{Statement of result} In this note we prove the
following: 

\begin{th}\label{main}
Let $\Gamma\subset PSL_2(\bfz)$ be a congruence subgroup, and $X_\Gamma$ the
corresponding modular curve. Let $D_\Gamma = [PSL_2(\bfz):\Gamma]$ and let
$d_\bfc(X_\Gamma)$ be the  $\bfc$-gonality of 
$X_\Gamma$. Then $${7\over 800} D_\Gamma \leq 
d_\bfc(X_\Gamma).$$ 

For $\Gamma = \Gamma_0(N)$ we have that $d_\bfc(X_{\Gamma_0(N)})$ is bounded
below by 
 ${7\over {800}}\cdot N$.
 Similarly, we
obtain  a  quadratic lower bound in $N$ for  $d_\bfc(X_{\Gamma_1(N)})$.
\end{th}  

\subsection{Remarks}

The proof, which was included in the author's thesis \cite{thesis}, follows
closely a suggestion of N. Elkies.  In the exposition here 
many details were added to the argument in \cite{thesis}. 

We utilize the work \cite{liyau} of P. Li
and S. T. Yau 
on conformal volumes, as well as the known bound on the leading nontrivial
eigenvalue of the non-euclidean Laplacian $\lambda_1\geq {{21}\over {100}}$
\cite{lrs}.  If Selberg's eigenvalue  
conjecture is true, the constant $7/800$ above may be replaced by $1/96$.

Since, by the Gauss - Bonnet formula, the genus $g(X_\Gamma)$ is bounded by
$D_\Gamma/12+1$ (indeed the difference is
$o(D_\Gamma)$), we may  rewrite the inequality above in the 
slightly weaker form $${{21}\over {200} } (g(X_\Gamma)-1) \leq
d_\bfc(X_\Gamma).$$  

For an analogous result about Shimura curves, see theorem \ref{shimura} below.

It should be noted (as was pointed out by P. Sarnak) that the gonality
has an {\em upper} bound of the same 
type. For the $\bfc$-gonality, by Brill-Noether theory \cite{kl} we
have $d_\bfc(X_\Gamma) \leq 1+\left[{{g+1}\over 2}\right]$. If, instead, 
one is interested in the gonality over the field of definition of $X_\Ga$, one
can use the canonical linear series to obtain the upper bound $2g-2$ if $g>1$,
and in the few cases where $g=1$ one can use the morphism to $X(1)$ and get the
upper bound 
$D_\Gamma$. 

\subsection{Acknowledgements} As mentioned above, I am indebted to Noam Elkies
for the main idea. The question was first brought to my attention in a letter
by S. Kamienny. The result first appeared in my thesis under the supervision of
Prof. J. Harris. Thanks are due to David Rohrlich and Glenn Stevens who set me 
straight on some details, and to Peter Sarnak for helpful suggestions.

\section{Setup and proof}
\subsection{Gonality}\label{gonal} Let $C$ be a smooth, projective, absolutely
irreducible 
algebraic curve over a field $K$.  Define the
$K$-{\bf gonality} $d_K(C)$ of $C$ to be the minimum degree of a finite
$K$-morphism $f:C\rar \bfp^1_K$. Clearly if $K\subset L$ then  
$d_K(C)\geq d_L(C\times_KL)$, and equality must hold whenever $K$ is
algebraically closed. 

\subsection{Congruence subgroups and modular curves}
By a {\bf congruence subgroup}  $\Gamma\subset PSL_2(\bfz)$ we mean that for
some $N$, $\Gamma$ 
contains the principal congruence subgroup $\Ga(N)$ of $2\times 2$ integer
matrices congruent to the identity modulo $N$.

Since $PSL_2(\bfr)$ acts on $\bfh= \{z = x+iy|y>0\}$ via fractional linear
transformations, we may let $Y_\Gamma = \GH$. It is well known that $Y_\Gamma$
may 
be compactified by adding finitely many points, called {\bf cusps}, to obtain a
compact Riemann surface $X_\Ga$, which we call the {\bf modular curve}
corresponding to $\Ga$. 

\subsection{The Poincar\'e metric}
The upper half plane $\bfh$ carries the Poincar\'e metric $ds^2 = {{dx^2 +
dy^2}\over {y^2}}$, which is  $PSL_2(\bfr)$ - invariant.  The
corresponding area form is given by ${{dx \, dy} \over {y^2}}$.
 Away from a finite set $T$
consisting the cusps and possibly some elliptic fixed points, the metric
descends to a Riemannian metric on $X_\Ga {\,\,^{_\setminus}\,\,} T$, of finite
area. 
We denote the area measure by $d\mu$. 

We will accordingly call a quadratic differential $ds^2$ a {\bf singular
metric} if it is a Riemannian metric away from finitely many points, and has
finite area. Thus the Poincar\'e metric gives rise to a singular metric on
$X_\Ga$. 

\subsection{The Laplacian} It is natural to consider the Hilbert space
$L_2(\GH)= L_2(X_\Gamma)$, where the $L_2$ pairing is taken with respect to the
Poincar\'e metric. The Laplace-Beltrami operator associated with the metric
$$ \Delta = -y^2({{\partial^2}\over {\partial x^2}} +{{\partial^2}\over
{\partial y^2}})$$ gives rise to a self adjoint unbounded operator on $
L_2(X_\Gamma)$, which is in fact positive semidefinite.

The kernel of $\Delta$ consists of the constant functions. In contrast with the
case of a genuine Riemannian metric on a compact manifold, the spectrum of
$\Delta$ is not discrete (see e.g. \cite{hejhal}, VI\S 9, VII\S 2, VIII\S 5).
The continuous spectrum is $\{\lambda\geq 1/4\} \subset \bfr$, and is fully
accounted for by an integral formula involving Eisenstein series $E(z,s)$ for
$Re(s) = 1/2$. The discrete part of the spectrum is given by $\lambda_0=0$
corresponding to the constants, and $0< \lambda_1 < \lambda_2...$ corresponding
to the so called {\bf cuspidal} eigenvectors.

\subsection{Selberg's conjecture} The question, what is $\lambda_1$ turns out
to 
be a fundamental one. Selberg \cite{selberg} has shown that $\lambda_1\geq
3/16$ and conjectured that $\lambda_1\geq 1/4$. Recently, Luo, Rudnick and
Sarnak \cite{lrs} showed that  $\lambda_1\geq 0.21$ (note that $3/16 < 0.21 <
1/4$). 

Since the continuous spectrum is known to be $\lambda\geq 1/4$, denote by
$\lambda_1' = \min(\lambda_1, 1/4)$. The value of $\lambda_1'$ has the
following characterization: 

Let $g$ be a nonzero continuous, piecewise differentiable function on $X_\Ga$
such that $\nabla g$ is square integrable with respect to $\mu$, and 
$\int_{X_\Ga} g d\mu=0$. Then (identifying $X_\Ga$ with $\GH$) we have

$$\int_{\GH} \left(\,\left({{\partial g}\over {\partial x}}\right)^2
+\left({{\partial 
g}\over {\partial  y}}\right)^2 \,\right) dx \,  dy \geq \lambda_1' \int_{\GH}
g^2 
{{dx \, dy} \over {y^2}}.$$
This is, in fact, the way Selberg originally stated his result.

\subsection{Conformal area} 
Let $C$ be a compact Riemann surface. Following
\cite{liyau}, we define the {\bf conformal area}, or the first conformal
volume $A_c(C)$ to be the infimum of $\int_C f^* d\mu_0$, where
$f:C\rar\bfp^1_\bfc$ runs over all nonconstant conformal mappings, and
where $d\mu_0$ is the $SO_3$-invariant area element on the Riemann sphere.
Using the conformal property of homotheties in $\bfp^1$, Li and Yau show
easily that $$A_c(C) \leq 4\pi\cdot d_\bfc(C).$$ 

On the other hand, given a Riemannian metric on $C$, let $A(C)$ be the area of
$C$. Using an elegant fixed point argument,  Li and Yau obtain (\cite{liyau},
Theorem 1)  $$\lambda_1 A(C) \leq 2A_c(C).$$  
Their proof works word for word in the case of our singular metric on
$X_\Ga$, once we replace $\lambda_1$ by $\lambda_1'$.  All that is
needed is, first, the characterization of $\lambda_1'$ discussed above, and
second, the
fact that differentiable functions on $X_\Gamma$ have a square-integrable
gradient. The latter follows since $\int_{X_\Gamma}|\nabla g|^2d\mu$ is
 invariant under conformal change of the metric, therefore it may be calculated
using a regular metric, and thus is finite.

\subsection{Conclusion of the proof} Since the Poincar\'e metric on $X_\Ga$ is
pulled back from $X_{PSL_2(\bfz)}=X(1)$,
we have $A(X_\Ga)=D_\Gamma \cdot A(X(1)) = D_\Gamma \cdot\pi/3$.
Combine this with the inequalities of Li and Yau, and obtain the first part of
the theorem.  Now note that $[PSL_2(\bfz):\Gamma_0(N)]$ is at least $N$, and
similarly $[PSL_2(\bfz):\Gamma_1(N)]$ is quadratic in $N$ (between $6(N/\pi)^2$
and $N^2$), and obtain the second part.
\qed

\subsection{An analogous result for Shimura curves} As was pointed out by
P. Sarnak, we have the following:

\begin{th}\label{shimura}
Let $D$ be an indefinite quaternion algebra over $\bfq$, and
let $G$ be the group of units of norm 1 in some order of $D$. Let
$\Gamma\subset G$ be a subgroup of 
finite index, and let $X_\Gamma=\Gamma \setminus \bfh$ be the corresponding
Shimura curve. Then $${{21}\over {200} } (g(X_\Gamma)-1) \leq
d_\bfc(X_\Gamma).$$   
\end{th}

{\bf Proof.} Since $X_\Gamma$ is compact, every automorphic form $g$ appearing
in $L^2(X_\Gamma)$ is cuspidal. It follows from the Jacquet - Langlands 
correspondence (see \cite{gelbart}, Theorem 10.1 and Remark 10.4) that unless
$g$ is the constant function, there exists a cuspidal automorphic form for some
congruence subgroup in $SL_2(\bfz)$ which has the same eigenvalue with respect
to the non-euclidean Laplacian. Therefore 
$\lambda_1 \geq 0.21$ holds for $X_\Gamma$. 
The results of Li and Yau give $\lambda_1 A(X_\Gamma) \leq 8\pi\cdot
d_\bfc(X_\Gamma)$, and the Gauss - Bonnet formula gives $4\pi(g(X_\Gamma))-1)
\leq A(X_\Gamma)$ (the difference coming from elliptic fixed points). Combining
the three inequalities we obtain the result. \qed

The author was informed that the results of \cite{lrs} were generalized by
Rudnick and Sarnak to
cuspidal automorphic forms on $GL_2$ over an arbitrary number field
$F$. Therefore Theorem \ref{shimura} holds for $D$ a quaternion algebra over a
totally real 
field, which is indefinite at exactly one infinite place.

\section{Applications and remarks}
\subsection{$\bfq$-gonality and rational torsion on elliptic curves} 
Let $C$ be a curve as in \ref{gonal}. Recall \cite{ah} that a point
$P\in C$ is called {\bf a point of degree $d$} if $[K(P):K]=d$. Suppose $C$ has
infinitely many points of degree $d$. By taking Galois orbits on the $d$-th
symmetric power of $C$ we have that $\sym^d(C)(K)$ is infinite. Let
$W_d(C)\subset Pic^d(C)$ be the image of $\sym^d(C)(K)$ by the Abel-Jacobi
map. In \cite{ah} 
it was noted that in this situation either $d_K(C) \leq d$, or $W_d(C)(K)$  
is infinite. 

Now assume $K$ is a number field. By a celebrated theorem of Faltings
\cite{fal}, if   $W_d(C)(K)$ is infinite then $W_d(C)\subset Pic^d(C)$ contains
a positive dimensional translate 
of an abelian variety, and the simple lemma 1 of \cite{ah} 
implies that $d_K(C)\leq 2d$ (\cite{thesis}, theorem 9). The latter
conclusion was also obtained by G. Frey in \cite{frey}.  

We now restrict attention to the case where $K=\bfq$ and $C = X_0(N)$. 
In
\cite{thesis}, Theorem 12, as well as in \cite{frey}, it was noted that a lower
bound on the $\bfq$-gonality, such as given by theorem \ref{main}, implies
that there exists a constant $m(d)$ (in fact, $m=230d$  will do), 
such that if $N> m(d)$ then $X_0(N)$ (and thus also $X_1(N)$) has finitely many
points of degree $d$. In section 1 of \cite{km}, Kamienny and Mazur showd that
this 
reduces the uniform boundedness conjecture on torsion points on elliptic curves
to bounding rational torsion of prime degree. The conjecture was finally proved
by L. Merel in \cite{merel}. 

It should be remarked that, since for this application one only needs a lower
bound on the $\bfq$-gonality of $X_0(N)$, one can use other methods, such as
Ogg's method \cite{ogg}. This is indeed the method used by Frey in \cite{frey},
although the bound obtained is not linear.
For points of low degree, one can use the main results of \cite{ah} with Ogg's
method to slightly improve the bound on $N$ (see \cite{hs} and \cite{thesis},
 2.5).  

For another arithmetic application of the lower bound on teh $\bfc$-gonality,
regarding pairs of elliptic curves with with isomorphic mod $N$
representations, see Frey \cite{frey2}.

\subsection{Torsion points: the function field case}
Recently, there has been renewed interest in the question of $\bfc$-gonality of
modular curves. In their paper \cite{ns}, K. V. Nguyen and M.-H. Saito used
algebraic techniques to give a lower bound on
the gonality. Although their bound is a bit weaker than ours, their methods are
of interest on their own right: they combine Ogg's method with a Castelnuovo
type bound. They pointed out that given any such bound, one obtains a
function field analogue of the strong uniform boundedness theorem about
torsion on elliptic curves, namely: given a non-isotrivial elliptic curve over
the function field of a complex curve  $B$, the size of the torsion
subgroup is bounded solely in terms of the gonality of $B$. This result is
strikingly analogous to a recent result of P. Pacelli (\cite{p}, Theorem 1.3):
assuming Lang's conjecture on rational curves on varieties of general type, 
the number of non-constant points on a  curve $C$ of genus $>1$
over the function field of $B$ is  bounded solely in terms of the genus of $C$
and the gonality of $B$.

\end{document}